\documentclass[12pt,a4paper]{article}

\usepackage{epsf,epsfig,amsfonts,amsgen,amsmath,amstext,amsbsy,amsopn,amsthm
}
\usepackage{amsmath,times,mathptmx}
\usepackage{amsfonts,amsthm,amssymb}
\usepackage{amsfonts}
\usepackage{graphics}
\usepackage{latexsym,bm}
\usepackage{amsfonts,amsthm,amssymb,bbding}
\usepackage{indentfirst}
\usepackage{graphicx,enumerate}
\usepackage{color}
\usepackage[colorlinks=true,anchorcolor=blue,filecolor=blue,linkcolor=blue,urlcolor=blue,citecolor=blue]{hyperref}
\usepackage{float}
\usepackage{tikz}
\newtheorem{thm}{Theorem}

\newtheorem{lemma}{Lemma}
\newtheorem{false statement}{False statement}

\theoremstyle{definition}

\newtheorem{claim}{Claim}

\newtheorem{conj}{Conjecture}
\newtheorem{corollary}{Corollary}

\newtheorem{problem}{Problem}

\begin{document}

\title{A spectral condition for the existence of cycles with consecutive odd lengths in non-bipartite
graphs
\footnote{Supported by the National Natural Science Foundation of China (Nos. 12011530064 and 11771141).}}
\author{{Zhiyuan Zhang, Yanhua Zhao}\thanks{Corresponding author. E-mail addresses: jsam0331@163.com (Z. Zhang); yhua030@163.com (Y. Zhao).}\\
{\footnotesize School of Mathematics, East China University of Science and Technology, Shanghai 200237, China}}
\date{}

\maketitle {\flushleft\large\bf Abstract:} A graph $G$ is called $H$-free, if it does not contain $H$ as a subgraph.
In 2010, Nikiforov proposed a  Brualdi-Solheid-Tur\'{a}n type problem: what is the maximum spectral radius of an $H$-free graph of order $n$?  In this paper, we consider the Brualdi-Solheid-Tur\'{a}n type problem for non-bipartite graphs. Let $K_{a, b}\bullet K_3$ denote the graph obtained by identifying a vertex of $K_{a,b}$ in the part of size $b$ and a vertex of $K_3$. We prove that  if $G$ is a non-bipartite graph of order $n$ satisfying $\rho(G)\geq \rho(K_{\lceil\frac{n-2}{2}\rceil, \lfloor\frac{n-2}{2}\rfloor}\bullet K_3)$, then $G$ contains all odd cycles  $C_{2l+1}$ for each integer $l\in[2,k]$ unless $G\cong K_{\lceil\frac{n-2}{2}\rceil, \lfloor\frac{n-2}{2}\rfloor}\bullet K_3$, provided that  $n$ is sufficiently large with respect to $k$. This resolves the problem posed by Guo, Lin and Zhao (2021).

\vspace{0.1cm}
\begin{flushleft}
\end{flushleft}

\section{Introduction}
All graphs considered here are simple and undirected and $k$ always denotes a positive integer.
Let $G$ be a graph with vertex set $V(G)$ ($n=|V(G)|$) and edge set $E(G)$ ($e(G)=|E(G)|$).
For any vertices $u,v \in V (G)$, we denote  by $N_G(u)$ the neighborhood of $u$ in $G$ and  $d_G(u)= |N_G(u)|$.
The \emph{adjacency matrix} of $G$ is an $n\times n$ matrix, denote by $A(G)$, whose $(i,j)$-entry is $1$ if $v_i$ is adjacent to $v_j$ and 0 otherwise. The \emph{spectral radius} of $G$, denoted by $\rho(G)$, is the largest modulus of the eigenvalues of $A(G)$.

A graph $G$ is called \emph{$H$-free}, if it does not contain a subgraph isomorphic to $H$.
For $n>k$, let $S_{n, k}$ be the join of a clique on $k$ vertices with  $n-k$ isolated vertices, and $S^+_{n,k}$ be the graph obtained from $S_{n,k}$ by adding one edge.
 In 2010, Nikiforov \cite{N10} investigated  how large $\rho(G)$ could  be when $G$
 contains no cycles and paths of specified order.
 He also raised a Brualdi-Solheid-Tur\'{a}n type problem: what is the maximum spectral radius of an $H$-free graph of order $n$? Besides, he conjectured the solution to the Brualdi-Solheid-Tur\'{a}n type problem for even cycles as follows.
\begin{conj}[\cite{N10}]\label{conj1}
For $k\geq2$ and sufficiently large $n$, $S^+_{n,k}$ is the unique graph attaining the maximum spectral radius among all graphs with no $C_{2k+2}$.
\end{conj}
Conjecture \ref{conj1} was solved for partial cases by Zhai and Wang \cite{ZW12},  Nikiforov \cite{N09}, Zhai and Lin \cite{ZL20} and completely solved by Cioab\u{a},  Desai and Tait \cite{CDT22+}.
 Also, much attention has been paid to Brualdi-Solheid-Tur\'{a}n type problem for various types of $H$, such as clique \cite{N07,WI}, path \cite{N10}, and complete bipartite graph \cite{BG,N07,NI4}, etc. For more relevant results, we refer the reader to \cite{ GH19,GN20,LN,LF22+,N11} and references therein.
In particular, Nikiforov \cite{N08} presented the following result, which gave a spectral condition for the existence of $C_{2k+1}$ in a graph.

\begin{thm}[\cite{N08}] \label{c2k+1}
Let $G$ be a graph of sufficiently large order $n$ with $\rho(G)> \sqrt{\lfloor\frac{n^2}{4}\rfloor}$. Then $G$
contains a cycle of length $l$ for every $l\leq \frac{n}{320}$.
\end{thm}
Ning and Peng \cite{NP20} slightly refined Theorem \ref{c2k+1} as $n\geq 160t.$ Very recently, by using  different methods,  Zhai and Lin \cite{ZLar} improved the result to $n\geq 7t$ by omitting the condition ``sufficiently large order $n$'', and Li and Ning \cite{LN} improved it to $n\geq 4t$.
Considering that the extremal graph for Theorem \ref{c2k+1} is a balanced complete bipartite graph, Guo, Lin and Zhao \cite{GLZ} posed a  natural  Brualdi-Solheid-Tur\'{a}n type problem for non-bipartite graphs as follows.

\begin{problem}[\cite{GLZ}]\label{pr1}
What is the maximum spectral radius of a $C_l$-free non-bipartite graph of order $n$?
\end{problem}
If $l$ is even, then the answer to Conjecture \ref{conj1} implies  Problem \ref{pr1}.
 If $l$ is odd, 
 Lin, Ning and Wu \cite{LNW21+} solved Problem \ref{pr1} for $l=3.$
 Let $a$, $b$ be two positive integers with $a+b=n-2$.
We denote by $K_{a, b}\bullet K_3$ the graph obtained by identifying a vertex of $K_{a,b}$ belonging to the part of size $b$ and a vertex of $K_3$. Recently, for a non-bipartite graph, Guo, Lin and Zhao \cite{GLZ} gave an answer to  Problem \ref{pr1} for $l=5.$  They further posed the following problem.
\begin{problem}[\cite{GLZ}]\label{prob1}
If $\rho(G) \ge \rho(K_{\lceil \frac{n-2}{2} \rceil,\lfloor \frac{n-2}{2} \rfloor} \bullet K_3 )$, $G \ncong K_{\lceil \frac{n-2}{2} \rceil,\lfloor \frac{n-2}{2} \rfloor} \bullet K_3 $ and $n$ is sufficiently large with respect to $k$, does $G$ always contain a $C_{2k+1}$?
\end{problem}

In this paper, we give a spectral condition for the existence of cycles with consecutive odd lengths in non-bipartite graphs.
\begin{thm}\label{thm1.3}
Let $k\geq 2$ and $G$ be a non-bipartite graph of order $n$, where $n$ is sufficiently large with respect to $k$.
If $\rho(G) \ge \rho(K_{\lceil \frac{n-2}{2} \rceil,\lfloor \frac{n-2}{2} \rfloor} \bullet K_3 )$, then $G$ contains all odd cycles  $C_{2l+1}$ for each integer $l \in [2, k]$ unless $G \cong K_{\lceil \frac{n-2}{2} \rceil,\lfloor \frac{n-2}{2} \rfloor} \bullet K_3 $.
\end{thm}
By Theorem \ref{thm1.3}, we obtain a positive answer to Problem \ref{prob1}.
\begin{corollary}\label{thm1.4}
Let $k\geq 2$ and $G$ be a non-bipartite graph of order $n$, where $n$ is sufficiently large with respect to $k$.
 If $\rho(G) \ge \rho(K_{\lceil \frac{n-2}{2} \rceil,\lfloor \frac{n-2}{2} \rfloor} \bullet K_3 )$, then $G$ contains a $C_{2k+1}$ unless $G \cong K_{\lceil \frac{n-2}{2} \rceil,\lfloor \frac{n-2}{2} \rfloor} \bullet K_3 $.
\end{corollary}

\section{Preliminaries}
Let $K_{1,n-2}^1$ be the tree obtained by appending a path $P_2$ to a pendant vertex of $K_{1,n-2}$, and let $ K_{n-1}^1$ be the graph obtained by appending a complete graph $K_{n-1}$ to a pendant vertex of a path $P_2$.
 Guo, Wang, Lou and Li \cite{GWL2021+}  gave a new upper bound on the spectral radius of a graph which sightly improves Hong's bound:

\begin{lemma}[\cite{GWL2021+}]\label{upperboundedge}
If $G$ is a connected simple graph with $n\geq10$ vertices and $m$ edges, then
$\rho(G) \le \sqrt{2m-n}$ unless $G\cong K_{1,n-1},~ K_{1,n-2}^1,~K_n,~ K_{n-1}^1$.
\end{lemma}
For any set $V_1\subseteq V(G)$, we denote by $G[V_1]$ the subgraph induced by $V_1$ in $G$,
and $G-V_1$ (or $G-G[V_1]$)  the subgraph  induced by $V(G)\setminus V_1$ in $G$.

\begin{lemma}[\cite{LN}]\label{deletingvertex0}
Let $G$ be a graph. For any $v \in V(G)$,  we have $\rho^2(G) \le \rho^2(G-\{v\})+2d_G(v)$.
\end{lemma}

\begin{lemma}[\cite{LN}]\label{Consecutive cyclesni}
Let $\varepsilon$ be real with $0<\varepsilon<\frac 1 4$. Then there exists an integer $N := N(\varepsilon)$
such that if $G$ is a graph on $n$ vertices with $n \ge N$ and $\rho(G)>
\sqrt{\lfloor \frac{n^2}{4}\rfloor}$, then $G$ contains all cycles $C_l$ with $l \in [3,(\frac{1}{4}-\varepsilon)n]$.
\end{lemma}

We write $H \subseteq G$ if graph $H$ is a subgraph of $G$.
Let $\delta(G)=\min\{d_G(u):u\in V(G)\}.$
 The following result is a reduced version of Theorem $1$ of \cite{NP2006}. Since we will frequently apply the result, we show it as a lemma here.
\begin{lemma}[\cite{NP2006}]\label{nonbipartite1/3}
	Let $G$ be a non-bipartite graph of sufficiently large order $n$, and  $\delta(G) \ge \frac{n}{3}$. Then $C_l \subseteq G$ for every integer $l\in[4, \delta(G)+1]$.
\end{lemma}

Recall that $K_{\lceil \frac{n-2}{2} \rceil,\lfloor \frac{n-2}{2} \rfloor} \bullet K_3 $ is the graph obtained by identifying a vertex of $K_{\lceil \frac{n-2}{2} \rceil,\lfloor \frac{n-2}{2} \rfloor}$ belonging to the part of size $\lfloor \frac{n-2}{2} \rfloor$ and a vertex of $K_3$.
\begin{lemma}[\cite{GLZ}]\label{nonbipartitec5}
	Let $G$ be a non-bipartite graph with order $n \ge 21$. If $\rho(G) \ge \rho(K_{\lceil \frac{n-2}{2} \rceil,\lfloor \frac{n-2}{2} \rfloor} \bullet K_3 )$, then $G$ contains a pentagon unless $G \cong K_{\lceil \frac{n-2}{2} \rceil,\lfloor \frac{n-2}{2} \rfloor} \bullet K_3 $.
\end{lemma}

\begin{lemma}[\cite{N08}]\label{subgraphH0}
	Let $0<\alpha<\frac{1}{4}$, $0<\beta<\frac{1}{2}$, $\frac{1}{2}-\frac{\alpha}{4}<\gamma\le 1$, $K\ge 0$. For a sufficiently large integer $n$, if $G$ is a graph of order $n$ with
\begin{equation*}
  \rho(G)>\gamma n-\frac{K}{n} ~and~ \delta(G)\le(\gamma-\alpha)n,
\end{equation*}
then there exists an induced subgraph $H$  of $G$ with $|V(H)|\ge(1-\beta)n$ and satisfying one of the following conditions:
\begin{align*}
   &(i) ~ \rho(H)>\gamma(1+\frac{\alpha\beta}{2})|V(H)|; ~(ii) ~ \rho(H)>\gamma|V(H)|~ and ~\delta(H)>(\gamma-\alpha)|V(H)|.
\end{align*}
\end{lemma}

\begin{lemma}\label{claim0.2'}
For every integer $i\in [1,4]$, let $z_i \in V(G_i)$. we define $G_{i+1}=G_i-\{z_i\}$  satisfying $d_{G_i}(z_{i})\le \frac{n}{8}$ , and for any vertex $z\in V(G_5)$, $d_{G_5}(z)> \frac{n}{8} $. Let $G=G_1$, $H=G_5$ and $|V(G)|=n$. If $\rho^2(G)>\frac{n^2-4n+3}{4} $ and  $\rho^2(H)\leq \frac{|V(H)|^2}{4}=\frac{(n-4)^2}{4}$, then for each $i\in \{1,2,3,4\},$
$$d_G(z_i)\ge \frac{n}{8}-8.$$
\end{lemma}
\begin{proof}
Assume to the contrary that there exists some vertex $z_j\in \{z_1, z_2, z_3, z_4\}$ such that
\begin{equation}\label{dzj}
d_G(z_j)< \frac{n}{8}-8.
\end{equation}
We first assert that
\begin{equation}\label{lowerbounddegree}
  \sum\limits_{i=1}^{4}d_G(z_i)\geq \frac{1}{2}(\rho^2(G)-\rho^2(H))>\frac{n}{2}-\frac{13}{8}.
\end{equation}
In fact, note that $d_{G_{i}}(z_{i})\le d_{G}(z_{i})$ for $i\in \{2,3,4\}.$ Then by Lemma \ref{deletingvertex0}, we have
\begin{align*}
   &\rho^2(G_2)\geq \rho^2(G)-2d_G(z_1),\\
  &\rho^2(G_3) \geq \rho^2(G_2)-2d_G(z_2), \\
  &\rho^2(G_4)\geq  \rho^2(G_3)-2d_G(z_3), \\
 &\rho^2(H)= \rho^2(G_5)\geq  \rho^2(G_4)-2d_{G}(z_4).
\end{align*}
Summing up all these inequalities, we have  $2\sum\limits_{i=1}^{4}d_G(z_i)\geq \rho^2(G)-\rho^2(H)$.
Thus,
$\sum\limits_{i=1}^{4}d_G(z_i)\geq \frac{1}{2}(\rho^2(G)-\rho^2(H))> \frac{n^2-4n+3}{8}-\frac{(n-4)^2}{8}=\frac{n}{2}-\frac{13}{8},$
as desired.
On the other hand, we have
 $$
  d_{G}(z_1)\leq \frac{n}{8} ~\textrm{and }~d_{G}(z_{i+1})-i\leq d_{G_{i+1}}(z_{i+1})\leq \frac{n}{8} ~\textrm{for} ~ i\in \{1,2,3\},$$
which follows that $ d_{G}(z_2)\leq \frac{n}{8}+1 ,~ d_{G}(z_3)\leq \frac{n}{8}+2,~  d_{G}(z_4)\leq \frac{n}{8}+3.$
Thus, the degree sum of any three vertices in $\{z_1, z_2, z_3, z_4\}$ does not exceed $\frac{3n}{8}+6$.
 Combining this  with (\ref{dzj}) and (\ref{lowerbounddegree}), we obtain that
 $$\frac{3n}{8}+6\geq \sum\limits_{i=1}^{4}d_G(z_i)-d_G(z_j)\geq \frac{n}{2}-\frac{13}{8}-(\frac{n}{8}-8)=\frac{3n}{8}+\frac{51}{8},$$
which deduces a contradiction, as desired.
\end{proof}

\section{\bf{Proof of Theorem \ref{thm1.3}}}

In this section, we shall give the proof of Theorem \ref{thm1.3}.
\vspace{3mm}

\noindent{\bf{Proof of Theorem \ref{thm1.3}}}.
First assume that $G$ is disconnected. Let $G_1,...,G_\eta$ be all components of $G$ for $\eta\ge 2$, and let
$G'$ be a connected graph obtained from $G$ by adding $\eta-1$ edges.
By the Rayleigh quotient and the Perron-Frobenius theorem,  $\rho(G')> \rho(G)\geq \rho(K_{\lceil \frac{n-2}{2} \rceil,\lfloor \frac{n-2}{2} \rfloor} \bullet K_3)$. Also note that $G'\ncong K_{\lceil \frac{n-2}{2} \rceil,\lfloor \frac{n-2}{2} \rfloor} \bullet K_3$ and $G$ contains a cycle of length $l$ if and only if $G'$ contains a cycle of length $l$ for any $l\geq 3$.  Then it remains to prove the result for $G'$. For this reason,  in what follows, we always assume  that $G$ is connected.

Note that $K_{\lceil \frac{n-2}{2} \rceil,\lfloor \frac{n-2}{2} \rfloor}$ is a proper subgraph of $K_{\lceil \frac{n-2}{2} \rceil,\lfloor \frac{n-2}{2} \rfloor} \bullet K_3$
and $\rho(K_{\lceil \frac{n-2}{2} \rceil,\lfloor \frac{n-2}{2} \rfloor})\geq \sqrt{\frac{n^2-4n+3}{4}}.$
Then
\begin{equation}\label{lowerboundG}
  \rho^2(G)\geq \rho^2(K_{\lceil \frac{n-2}{2} \rceil,\lfloor \frac{n-2}{2} \rfloor} \bullet K_3 )>\rho^2(K_{\lceil \frac{n-2}{2} \rceil,\lfloor \frac{n-2}{2} \rfloor})\geq
   \frac{n^2-4n+3}{4}.
\end{equation}
It is clear that $\min \{ \rho(K_n), \rho(K_{n-1}^1)\}>n-2> \rho(K_{\lceil \frac{n-2}{2} \rceil,\lfloor \frac{n-2}{2} \rfloor} \bullet K_3)$, and graphs $K_n$ and  $K_{n-1}^1$ contain all odd cycles $C_{2l+1}$ for each integer $l\in [2,k]$, where $n$ is sufficiently large with respect to $k$.
Then we consider that $G \notin  \{K_{1,n-1},~ K_{1,n-2}^1\}$ in the following. Furthermore, since $G$ is non-bipartite graph, we have $G\notin \{K_{1,n-1},~ K_{1,n-2}^1\}.$
Combining this  with Lemma \ref{upperboundedge} and (\ref{lowerboundG}), we obtain $2e(G) > \frac{n^2+3}{4}$, which implies that the average degree of $G$ is $\overline{d}(G)=\frac{2e(G)}{n} >\frac{n}{4}$.
Let $H$ be an induced subgraph of $G$  defined by  a sequence of graphs $G_1, G_2, \cdots, G_{t+1}$ such that:
\begin{enumerate}[(1)]
    \item $G_1=G$, $G_{t+1}=H$;
  \item for every integer $j\in [1,t]$ (if $t\geq 1$), there exists one vertex $v_{i_j}\in V(G_j)$ such that $d_{G_j}(v_{i_j})\le \frac{n}{8}$ and $G_{j+1}=G_j-\{v_{i_j}\}$;
  \item for every vertex $v\in V(G_{t+1})$ (if any), $d_{G_{t+1}}(v)>\frac{n}{8}$.
\end{enumerate}
By the construction, $|V(G)\setminus V(H)|=t$, and if $t\ge 1$ then $d_{G_j}(v_{i_{j}})\le \frac{n}{8}$ for $1\leq j\leq t$. We assert that $H$ is not a null graph. Otherwise, assume that $H$ is  a null graph, i.e., $t=n$. Then $$\overline{d}(G)=\frac{2(d_G(v_{i_1})+d_{G_2}(v_{i_2})+\cdots+d_{G_t}(v_{i_t}))}{n}\le \frac{2\times t\times \frac{n}{8}}{n}=\frac{n}{4},$$
contrary to $\overline{d}(G)>\frac{n}{4}$. Therefore, $G$ contains an induced subgraph $H$  with $\delta(H)>\frac{n}{8}$.

 Let $V_1=V(G)\setminus V(H)$ and $|V(H)|=h$. Then $|V_1|=t$ and $h+t=n.$ We assert that
\begin{equation}\label{lowerboundH}
\rho^2(H)>\frac {n^2-(4+t)n+3}{4}= \frac {(h+t)(h-4)+3}{4}.
\end{equation}
In fact, if $t=0$, then $H=G$, and the inequality holds by (\ref{lowerboundG}). If $t\ge 1$, since  $d_{G_j}(v_{i_{j}})\le \frac{n}{8}$ for every integer $j\in [1,t]$,
by  Lemma \ref{deletingvertex0}, we have
\begin{align*}
   & \rho^2(G_2) \geq \rho^2(G) - 2\cdot\frac{ n}{8},\\
  & \rho^2(G_3) \geq \rho^2(G_2)-2\cdot\frac{n}{8}, \\
   &~~~~~~~~~~~~~~~~~~~~~~~~~~~~~\vdots \\
 &\rho^2(G_t) \geq \rho^2(G_{t-1})-2\cdot\frac{n}{8}, \\
 & \rho^2(G-V_1)= \rho^2(H) \geq \rho^2(G_t)-2\cdot\frac{n}{8}.
\end{align*}
Then it follows from  (\ref{lowerboundG})  that
\begin{equation*}
 \rho^2(H) \ge \rho^2(G) - 2t\cdot\frac{n}{8} >\frac {n^2-(4+t)n+3}{4},
\end{equation*}
as desired.

Note that $|V_1|=t\ge 0.$ Then we shall distinguish our argument into $t\ge 5$ and $t\le 4$ these two situations.
One can verify that $h$ is sufficiently large with respect to $k$ since $h\geq \delta(H)+1 >\frac n 8 $.

For $t\geq 5$, by (\ref{lowerboundH}), we deduce that
$$\rho^2(H) > \frac {(h+t)(h-4)+3}{4}\geq \frac {(h+5)(h-4)+3}{4}=\frac {h^2+h-17}{4}\ge \frac {h^2}{4}.$$
  By Lemma \ref{Consecutive cyclesni}, $H$ contains all cycles $C_l$ for each  integer ${\mathit{l}} \in [3,(\frac{1}{4}-\epsilon)h]$, where $0<\epsilon <\frac{1}{4}$.
 Then we obtain that $H$ contains all odd cycles $C_{2l+1}$ for each integer $l\in [2,k],$ as required.

For $t\leq 4$, i.e., $h=n-t\geq n-4$.
Let\begin{center}
  $0<\theta\leq \frac{1}{321}$, $\alpha=\frac{1}{10}+\theta$, $\beta=40\theta$, $\gamma=\frac{1}{2}-\theta$ and $K=0$.
\end{center}
Combining this with (\ref{lowerboundH}), we have
\begin{equation}\label{lowerboundH2}
  \rho^2(H)>\frac {(h+t)(h-4)+3}{4}\geq \frac {h(h-4)+3}{4}>(\gamma h)^2=(\frac{1}{2}-\theta)^2h^2
\end{equation}
 for sufficiently large $h$.
 Furthermore, if $\rho^2(H)> \frac{h^2}{4}$, then by Lemma \ref{Consecutive cyclesni} again, we obtain that $H$ contains all odd cycles $C_{2l+1}$ for each integer $l\in [2,k],$ as required.
 In what follows, we suppose that
 \begin{equation}\label{upperbound}
   \rho^2(H)\le \frac{h^2}{4}.
 \end{equation}
Recall that $\delta(H)>\frac{n}{8}$.  Now we shall divide the proof into  two cases as follows. 

\vspace{3mm}
\noindent{\bf{$\underline{\mbox{Case 1. $\delta(H)>(\gamma-\alpha)h=(\frac{2}{5}-2\theta)h.$}}$}}
\vspace{3mm}

 Note that $\delta(H)>(\frac{2}{5}-2\theta)h>\frac{h}{3}.$ If $H$ is a non-bipartite graph, then  $H$ contains all cycles $C_{l}$ for each  integer $l\in [4,\delta(H)+1]$ by Lemma \ref{nonbipartite1/3}. Since $h$ is sufficiently large with respect to $k$, we obtain that $H$ contains all odd cycles $C_{2l+1}$ for each integer $l\in [2,k],$ as required.

Now we suppose that $H$ is a bipartite graph. Denote by $V_{1H}$ and $V_{2H}$ the two parts of $H$.
We first give  the following claim.
\begin{claim}\label{claim0.1}
 For any  vertex $u, v\in V_{1H}$ (or $ u, v\in V_{2H}$), $|N_H(u)\cap N_H(v)|>k+5.$
\end{claim}
\renewcommand\proofname{\bf Proof}
\begin{proof}
Note that $\delta(H)>(\frac{2}{5}-2\theta)h$ and $H$ is a bipartite graph. Then $\min\{|V_{1H}|, |V_{2H}|\}> (\frac{2}{5}-2\theta)h$ , which leads to $\max\{|V_{1H}|, |V_{2H}|\}\leq (\frac{3}{5}+2\theta)h$ since $|V_{1H}|+|V_{2H}|=h$. Without loss of generality, we suppose that $u, v\in V_{1H}$.
Then,
$$d_H(u)+d_H(v)-|V_{2H}|\geq (\frac{4}{5}-4\theta)h-(\frac{3}{5}+2\theta)h=(\frac{1}{5}-6\theta)h> k+5 $$
for sufficiently large $h$, which implies $|N_H(u)\cap N_H(v)|>k+5,$ as required.
\end{proof}

  Denote by $P_{uv}$  the path connecting the vertex $u$ to vertex $v$.
\begin{claim}\label{claim0.2}
Let $S\subseteq V(H)$ and $|S|\le 5.$ For any vertex $u_1\in V_{1H}$ and $ v_1\in V_{2H}$,
there exists a path $P_{u_1v_1}$ of order $2l$ for each integer $ l\in [2, k]$ in $H-S$.
 \end{claim}
\renewcommand\proofname{\bf Proof}
\begin{proof}
Since $\delta(H)> (\frac{2}{5}-2\theta)h$ and $h$ is sufficiently large, there exists a vertex $v_{2}\in N_{H-S}(u_1)\setminus \{v_1\}$  (\emph{See Figure \ref{Figure 1}.}). By Claim \ref{claim0.1}, we can find a vertex $u_2\in N_{H-S}(v_1)\cap (N_{H-S}(v_{2})\setminus \{u_1\})$.  By the same reason, there exists a vertex $v_{3}\in N_{H-S}(u_{2})\setminus \{v_1, v_{2}\}$ and  a vertex $u _{3}\in N_{H-S}(v_{3})\cap (N_{H-S}(v_{2})\setminus \{u_1, u_{2}\})$.   Repeating the  procedure, we find  a vertex $v_{k}\in N_{H-S}(u_{k-1})\setminus \{v_1, \cdots, v_{k-1}\}$ and a vertex $u_{k}\in N_{H-S}(v_{k-1})\cap (N_{H-S}(v_{k})\setminus \{u_1,\cdots, u_{k-1}\})$. Thus, ${H-S}$ contains a path $P_{u_1v_1}:=u_1v_2u_3v_4\cdots u_{l-1}v_lu_lv_{l-1}\cdots u_2v_1$ of order $2l$  for each even integer $ l\in[2, k]$,
or a path $P_{u_1v_1}:=u_1v_2u_3v_4\cdots v_{l-1}u_lv_lu_{l-1}\cdots u_2v_1$ of order $2l$  for each odd integer $ l\in[3, k]$,
as required.
\end{proof}
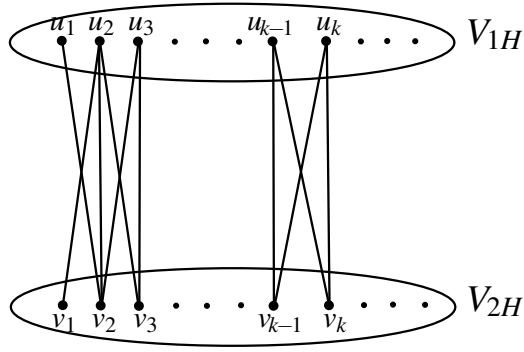
\begin{figure}[hbtp]
\centering
\begin{tikzpicture}[x=1.00mm, y=1.00mm, inner xsep=0pt, inner ysep=0pt, outer xsep=0pt, outer ysep=0pt,,scale=0.5]
\path[line width=0mm] (39.04,16.93) rectangle +(136.63,94.73);
\definecolor{L}{rgb}{0,0,0}
\path[line width=0.30mm, draw=L] (99.50,99.37) ellipse (58.45mm and 10.29mm);
\definecolor{F}{rgb}{0,0,0}
\path[line width=0.30mm, draw=L, fill=F] (54.69,99.84) circle (1.00mm);
\path[line width=0.30mm, draw=L, fill=F] (64.67,99.75) circle (1.00mm);
\path[line width=0.30mm, draw=L, fill=F] (74.73,99.66) circle (1.00mm);
\path[line width=0.30mm, draw=L, fill=F] (84.43,99.59) circle (0.50mm);
\path[line width=0.30mm, draw=L, fill=F] (92.61,99.68) circle (0.50mm);
\path[line width=0.30mm, draw=L, fill=F] (101.33,99.57) circle (0.50mm);
\path[line width=0.30mm, draw=L, fill=F] (110.13,99.68) circle (1.00mm);
\path[line width=0.30mm, draw=L, fill=F] (124.21,99.76) circle (1.00mm);
\path[line width=0.30mm, draw=L, fill=F] (133.36,99.66) circle (0.50mm);
\path[line width=0.30mm, draw=L, fill=F] (140.29,99.77) circle (0.50mm);
\path[line width=0.30mm, draw=L, fill=F] (147.66,99.72) circle (0.50mm);
\path[line width=0.30mm, draw=L] (99.75,29.23) ellipse (58.45mm and 10.29mm);
\path[line width=0.30mm, draw=L, fill=F] (54.95,29.70) circle (1.00mm);
\path[line width=0.30mm, draw=L, fill=F] (64.92,29.61) circle (1.00mm);
\path[line width=0.30mm, draw=L, fill=F] (74.98,29.52) circle (1.00mm);
\path[line width=0.30mm, draw=L, fill=F] (84.77,29.58) circle (0.50mm);
\path[line width=0.30mm, draw=L, fill=F] (92.87,29.54) circle (0.50mm);
\path[line width=0.30mm, draw=L, fill=F] (101.85,29.70) circle (0.50mm);
\path[line width=0.30mm, draw=L, fill=F] (110.39,29.54) circle (1.00mm);
\path[line width=0.30mm, draw=L, fill=F] (124.97,29.68) circle (1.00mm);
\path[line width=0.30mm, draw=L, fill=F] (134.19,29.73) circle (0.50mm);
\path[line width=0.30mm, draw=L, fill=F] (141.81,29.86) circle (0.50mm);
\path[line width=0.30mm, draw=L, fill=F] (150.25,29.81) circle (0.50mm);
\path[line width=0.30mm, draw=L] (54.78,99.89) -- (65.02,29.63);
\path[line width=0.30mm, draw=L] (64.66,99.85) -- (54.96,29.58);
\path[line width=0.30mm, draw=L] (64.48,99.49) -- (65.20,29.58);
\path[line width=0.30mm, draw=L] (64.66,100.03) -- (74.91,29.22);
\path[line width=0.30mm, draw=L] (75.27,99.31) -- (75.09,29.76);
\path[line width=0.30mm, draw=L] (74.91,99.67) -- (65.20,29.58);
\path[line width=0.30mm, draw=L] (110.31,99.67) -- (110.49,29.58);
\draw(51.56,102.67) node[anchor=base west]{\fontsize{11.23}{17.07}\selectfont$u_{1}$};
\draw(61.60,102.67) node[anchor=base west]{\fontsize{11.23}{17.07}\selectfont $u_{2}$};
\draw(72.15,102.67) node[anchor=base west]{\fontsize{11.23}{17.07}\selectfont $u_{3}$};
\draw(103.45,102.67) node[anchor=base west]{\fontsize{11.23}{17.07}\selectfont \selectfont $u_{\!k\!-\!1\!}$};
\draw(121.58,102.67) node[anchor=base west]{\fontsize{11.23}{17.07}\selectfont $u_{k}$};
\path[line width=0.30mm, draw=L] (110.38,99.65) -- (125.12,29.20);
\path[line width=0.30mm, draw=L] (125.12,29.56) -- (124.40,100.01);
\path[line width=0.30mm, draw=L] (124.40,100.01) -- (110.74,29.92);
\draw(52.07,23.91) node[anchor=base west]{\fontsize{11.23}{17.07}\selectfont $v_{1}$};
\draw(62.62,23.91) node[anchor=base west]{\fontsize{11.23}{17.07}\selectfont $v_{2}$};
\draw(73.53,23.91) node[anchor=base west]{\fontsize{11.23}{17.07}\selectfont $v_{3}$};
\draw(106.28,23.91) node[anchor=base west]{\fontsize{11.23}{17.07}\selectfont $v_{\!k\!-\!1\!}$};
\draw(123.47,23.91) node[anchor=base west]{\fontsize{11.23}{17.07}\selectfont $v_{k}$};
\draw(161.36,99.10) node[anchor=base west]{\fontsize{14.23}{17.07}\selectfont $V_{1H}$};
\draw(161.45,28.30) node[anchor=base west]{\fontsize{14.23}{17.07}\selectfont $V_{2H}$};
\end{tikzpicture}%
\caption{The path $P_{u_{1}v_{1}}$ in $H$.}\label{Figure 1}
\end{figure}

Furthermore,  we give the following result.
\begin{claim}\label{findcycle}
For any vertex $u_1\in V_{1H}$ and $ v_1\in V_{2H}$, let $P_{u_1v_1}\subseteq G.$
 If $|V(P_{u_1v_1})|=5$ , then $G$ contains all odd cycles $C_{2l+3}$ for each integer $l\in[2,k-1]$.
 \end{claim}
\renewcommand\proofname{\bf Proof}
\begin{proof}
Denote by $P_{u_1v_1}:=u_1w_1w_2w_3v_1$ and $S=V(H)\cap (V(P_{u_1v_1})\setminus \{u_1, v_1\})$. Then by Claim \ref{claim0.2}, there exists a path $P'_{u_1v_1}$ of order $2l$ for each integer $ l\in [2, k-1]$ in $H-S$.
Thus, $G$ contains all odd cycles $C_{2l+3}:=u_1w_1w_2w_3v_1P'_{v_1u_1}$ for each integer $l\in[2, k-1]$, as desired.
\end{proof}

Since $\rho(G)\geq \rho(K_{\lceil \frac{n-2}{2} \rceil,\lfloor \frac{n-2}{2} \rfloor} \bullet K_3 )$ and $G\ncong K_{\lceil \frac{n-2}{2} \rceil,\lfloor \frac{n-2}{2} \rfloor} \bullet K_3 $, by Lemma \ref{nonbipartitec5}, $G$ contains  a cycle $C_5:=z_1z_2z_3z_4z_5z_1$  as a subgraph.
In the following, we shall prove that $G$ also contains all odd cycles  $C_{2l+3}$ for each integer $l\in[2,k-1]$ to complete the proof of Case 1.

Note that $H$ is an induced  bipartite subgraph of $G$. Then $C_5$ is not a subgraph of $H$, and so $|V(C_5)\cap V(H)|\leq 4.$ Since $t=|V_1|=|V(G-H)|=n-h\leq 4$, we deduce that $|V(C_5)\cap V(H)|\ge 1.$ Therefore,  $1\leq |V(C_5)\cap V(H)|\leq 4.$
Now, we shall distinguish three situations as follows.
\begin{equation*}
  (i) ~ 3\leq |V(C_5)\cap V(H)|\leq 4, ~(ii)~|V(C_5)\cap V(H)|=2, ~(iii)~|V(C_5)\cap V(H)|=1.
\end{equation*}

For (i), there exist two vertices $u,v\in \{z_i| 1\le i\le 5\}$ such that $uv\in E(C_5)$ and  $u,v$ belong  to $V_{1H}$ and $V_{2H}$, respectively. Then we can find a path $P_{uv}$ of order $5$ in $G$, whose vertices belong to $V(C_5)$, and $u\in V_{1H}$ and $v \in V_{2H}$.
By Claim \ref{findcycle}, the result follows.

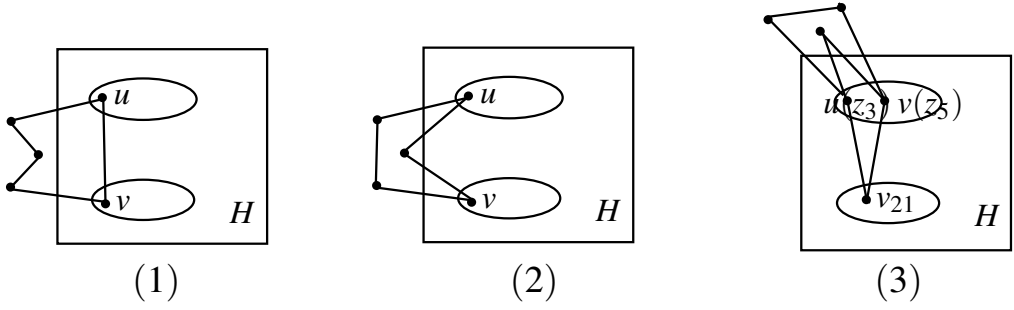
\begin{figure}[h]
\centering
\begin{tikzpicture}[x=1.30mm, y=1.30mm, inner xsep=0pt, inner ysep=0pt, outer xsep=0pt, outer ysep=0pt,scale=0.45]
\path[line width=0mm] (40.53,41.25) rectangle +(229.92,60.08);
\definecolor{L}{rgb}{0,0,0}
\definecolor{F}{rgb}{0,0,0}
\path[line width=0.30mm, draw=L, fill=F] (43.68,72.47) circle (1.00mm);
\path[line width=0.30mm, draw=L, fill=F] (43.53,57.57) circle (1.00mm);
\path[line width=0.30mm, draw=L, fill=F] (49.79,64.92) circle (1.00mm);
\path[line width=0.30mm, draw=L, fill=F] (64.17,77.86) circle (1.00mm);
\path[line width=0.30mm, draw=L, fill=F] (64.89,53.78) circle (1.00mm);
\path[line width=0.30mm, draw=L] (73.34,77.50) ellipse (15.63mm and 6.11mm);
\path[line width=0.30mm, draw=L] (73.34,54.68) ellipse (14.92mm and 5.93mm);
\path[line width=0.30mm, draw=L] (64.17,77.50) -- (42.96,72.47);
\path[line width=0.30mm, draw=L] (50.02,64.68) -- (43.53,72.18);
\path[line width=0.30mm, draw=L] (50.27,64.68) -- (43.28,57.18);
\path[line width=0.30mm, draw=L] (42.90,57.18) -- (64.63,54.26);
\path[line width=0.30mm, draw=L] (64.37,78.15) -- (64.88,54.01);
\path[line width=0.30mm, draw=L] (54.08,88.83) [rotate around={270:(54.08,88.83)}] rectangle +(44.09,47.14);
\draw(66.79,76.37) node[anchor=base west]{\fontsize{12.23}{17.07}$ u$};
\draw(67.17,52.48) node[anchor=base west]{\fontsize{12.23}{17.07}$v$};
\draw(92.84,49.43) node[anchor=base west]{\fontsize{12.23}{17.07}$H$};
\draw(70.85,33.90) node[anchor=base west]{\fontsize{14.23}{17.07}\selectfont $(1)$};
\path[line width=0.30mm, draw=L, fill=F] (125.99,72.83) circle (1.00mm);
\path[line width=0.30mm, draw=L, fill=F] (125.84,57.93) circle (1.00mm);
\path[line width=0.30mm, draw=L, fill=F] (132.10,65.28) circle (1.00mm);
\path[line width=0.30mm, draw=L, fill=F] (146.47,78.22) circle (1.00mm);
\path[line width=0.30mm, draw=L, fill=F] (147.19,54.14) circle (1.00mm);
\path[line width=0.30mm, draw=L] (155.64,77.86) ellipse (15.63mm and 6.11mm);
\path[line width=0.30mm, draw=L] (155.64,55.04) ellipse (14.92mm and 5.93mm);
\path[line width=0.30mm, draw=L] (146.47,77.86) -- (125.27,72.83);
\path[line width=0.30mm, draw=L] (125.20,57.54) -- (146.93,54.62);
\path[line width=0.30mm, draw=L] (136.39,89.18) [rotate around={270:(136.39,89.18)}] rectangle +(44.09,47.14);
\draw(149.09,76.73) node[anchor=base west]{\fontsize{12.23}{12.07}$u$};
\draw(149.47,52.84) node[anchor=base west]{\fontsize{12.23}{12.07}$ v$};
\draw(175.14,49.79) node[anchor=base west]{\fontsize{12.23}{12.07}$H$};
\draw(155.80,33.90) node[anchor=base west]{\fontsize{14.23}{17.07}\selectfont $(2)$};
\path[line width=0.30mm, draw=L] (221.31,87.34) [rotate around={270:(221.31,87.34)}] rectangle +(44.09,47.14);
\path[line width=0.30mm, draw=L, fill=F] (213.84,95.52) circle (1.00mm);
\path[line width=0.30mm, draw=L, fill=F] (230.33,98.27) circle (1.00mm);
\path[line width=0.30mm, draw=L, fill=F] (240.05,77.13) circle (1.00mm);
\path[line width=0.30mm, draw=L, fill=F] (231.65,77.14) circle (1.00mm);
\path[line width=0.30mm, draw=L, fill=F] (235.93,54.72) circle (1.00mm);
\path[line width=0.30mm, draw=L] (240.82,76.79) ellipse (15.63mm and 6.11mm);
\path[line width=0.30mm, draw=L] (240.82,53.96) ellipse (14.92mm and 5.93mm);
\draw(226.27,74.27) node[anchor=base west]{\fontsize{12.23}{17.07}$u(z_3)$};
\draw(241.15,74.27) node[anchor=base west]{\fontsize{12.23}{17.07} $v(z_5)$};
\draw(260.32,48.71) node[anchor=base west]{\fontsize{12.23}{17.07}$H$};
\draw(237.76,33.90) node[anchor=base west]{\fontsize{14.23}{17.07}\selectfont $(3)$};
\path[line width=0.30mm, draw=L] (146.24,77.90) -- (132.51,66.08);
\path[line width=0.30mm, draw=L] (131.75,65.32) -- (147.26,54.77);
\path[line width=0.30mm, draw=L] (126.03,72.56) -- (125.65,57.69);
\path[line width=0.30mm, draw=L, fill=F] (225.56,92.89) circle (1.00mm);
\path[line width=0.30mm, draw=L] (214.38,95.27) -- (230.77,78.37);
\path[line width=0.30mm, draw=L] (225.82,93.36) -- (231.28,77.86);
\path[line width=0.30mm, draw=L] (225.82,92.60) -- (239.92,77.35);
\path[line width=0.30mm, draw=L] (240.05,77.22) -- (229.76,99.33);
\path[line width=0.30mm, draw=L] (230.27,98.32) -- (213.75,96.16);
\path[line width=0.30mm, draw=L] (236.62,54.41) -- (236.49,53.91);
\path[line width=0.30mm, draw=L] (231.66,77.29) -- (235.98,54.92);
\path[line width=0.30mm, draw=L] (240.05,77.03) -- (235.98,54.16);
\draw(238.14,53.02) node[anchor=base west]{\fontsize{12.23}{10.07}$v_{21}$};
\end{tikzpicture}%
\caption{ The situations of Subcase 2.1(ii)}\label{fig00}
\end{figure}

For (ii), let $V(C_5)\cap V(H)=\{u,v\}$. Clearly, $\{u,v\}\subseteq \{z_i| 1\le i\le 5\}$.
If $uv\in E(C_5)$ for $u\in V_{1H}$, $v\in V_{2H}$ (\emph{See Figure \ref{fig00}.(1)}),  then by Claim \ref{findcycle},  $G$ contains all odd cycles $C_{2l+3}$ for each integer $l \in [2,k-1]$, as required.
If $uv\notin E(C_5)$ for $u\in V_{1H}$, $v\in V_{2H}$ (\emph{See Figure \ref{fig00}.(2)}), then by Claim \ref{claim0.2},  $G$ contains all odd cycles $C_{2l+1}$ for each integer $l \in [2,k]$, as required.
If  $u, v\in V_{1H}$ ( the situation $u, v\in V_{2H}$ can be proof by using the same analysis) (\emph{See Figure \ref{fig00}.(3)}), then there exists a vertex $v_{21}\in V_{2H}$ such that $v_{21}\in N_H(u)\cap N_H(v)$ by Claim \ref{claim0.1}.
 Without loss of generality, let $u=z_3$ and $v=z_5$. Then we can find a path $P_{z_3v_{21}}:=z_3z_2z_1z_5v_{21}$ of order $5$ in $G$ satisfying $z_3\in V_{1H}$ and $v_{21} \in V_{2H}$.
By Claim \ref{findcycle}, we complete the proof of (ii).

 For (iii),  we suppose $V(C_5)\cap V(H)=V(C_5)\cap V_{1H}=\{z_5\}$, and so $\{z_1,z_2,z_3,z_4\}\subseteq V(G-H)$.  Recall  $|V(G-H)|=|V_1|=t\le 4.$
Then $\{z_1,z_2,z_3,z_4\}=V(G-H).$
 By (\ref{upperbound}) and Lemma \ref{claim0.2'}, we deduce that for each $i\in \{1, 2, 3, 4\}$,
$$d_H(z_i)\geq d_G(z_i)-|V(G-H)|  \geq \frac{n}{8}-8-4 >2$$
for sufficiently large $n.$
Hence, there exists a vertex $v_{11}\in N_H(z_3)\setminus \{z_5\}$.
If  $v_{11}\in  V_{1H}$, then by Claim \ref{claim0.1}, there exists a vertex $v_{21}\in N_H(z_5)\cap N_H(v_{11}).$
 Thus we can find a path $P_{z_5v_{21}}:=z_5z_4z_3v_{11}v_{21}$ of order $5$ in $G$ satisfying $z_5\in V_{1H}$ and $v_{21} \in V_{2H}$.
By Claim \ref{findcycle}, the result follows.
If $v_{11}\in V_{2H}$,  then we can find a path $P_{z_5v_{11}}:=z_5z_1z_2z_3v_{11}$ of order $5$ in $G$ satisfying $z_5\in V_{1H}$ and $v_{11} \in V_{2H}$. By Claim \ref{findcycle}, we complete the proof of (iii).

\vspace{3mm}
\noindent{\bf{$\underline{\mbox{Case 2. $\frac{n}{8}<\delta(H)\leq(\gamma-\alpha)h=(\frac{2}{5}-2\theta)h.$}}$}}
\vspace{3mm}

Recall that
  $0<\theta\leq \frac{1}{321}$, $\alpha=\frac{1}{10}+\theta$, $\beta=40\theta$, $\gamma=\frac{1}{2}-\theta$ and $K=0$.
By (\ref{lowerboundH2})  and  Lemma \ref{subgraphH0},
we can find an induced subgraph $H_0$ of $H$ with $|V(H_0)|\geq(1-\beta)h=(1-40\theta)h\geq (1-40\theta)(n-4)>\frac{280}{321}n$ and satisfying one of the following conditions:
\begin{description}
  \item[(\textbf{I})] $\rho(H_0)>\gamma(1+\frac{\alpha\beta}{2})|V(H_0)|=(\frac{1}{2}+8\theta^2-20\theta^3)|V(H_0)|$,
  \item[(\textbf{II})]  $\rho(H_0)>\gamma|V(H_0)|=(\frac{1}{2}-\theta)|V(H_0)|$ and $\delta(H_0)>(\gamma-\alpha)|V(H_0)|=(\frac{2}{5}-2\theta)|V(H_0)|$.

\end{description}
 Note that $ |V(H_0)|$ is sufficiently large with respect to $k$ since $|V(H_0)|>\frac{280}{321}n$.
For \textbf{(I)}, we see that $\rho(H_0)>\frac{|V(H_0)|}{2}$, which follows that $H_0$ contains all cycles $C_{l}$ for each integer $ l\in[3,(\frac{1}{4}-\epsilon)|V(H_0)|]$  by Lemma \ref{Consecutive cyclesni}, where $0<\epsilon<\frac{1}{4}$. Thus, $H_0$ contains all odd cycles $C_{2l+1}$ for each integer $l\in [2,k] $, as required.

For \textbf{(II)}, the proof goes like in part Case 1, but needs more care.
Note that $\delta(H_0)>(\frac{2}{5}-2\theta)|V(H_0)|>\frac{|V(H_0)|}{3}.$
If $H_0$ is a non-bipartite graph,
then $H_0$ contains all cycles $C_{l}$ for each integer $l\in [4,\delta(H_0)+1]$ by Lemma \ref{nonbipartite1/3}. Thus, $H_0$ contains all odd cycles $C_{2l+1}$ for each integer $l\in [2,k] $, as required.

Now we focus on the situation that $H_0$ is a bipartite graph. Denote by $V_{1H_0}$ and $V_{2H_0}$ the two parts of $H_0$ (\emph{See Figure \ref{fig0}.}).

\begin{figure}[h]
\centering
\begin{tikzpicture}[x=1.00mm, y=1.00mm, inner xsep=0pt, inner ysep=0pt, outer xsep=0pt, outer ysep=0pt,scale=0.5]
\path[line width=0mm] (28.54,53.35) rectangle +(79.79,49.85);
\definecolor{L}{rgb}{0,0,0}
\path[line width=0mm] (37.84,38.41) rectangle +(84.44,63.23);
			\definecolor{L}{rgb}{0,0,0}
			\path[line width=0.30mm, draw=L] (39.84,99.64) [rotate around={270:(39.84,99.64)}] rectangle +(59.23,80.44);
			\path[line width=0.30mm, draw=L] (44.99,94.29) [rotate around={270:(44.99,94.29)}] rectangle +(44.33,64.77);
			\path[line width=0.30mm, draw=L] (75.18,74.99) [rotate around={360:(75.18,74.99)}] ellipse (25.03mm and 14.90mm);
			\path[line width=0.30mm, draw=L] (49.96,75.37) -- (100.21,75.37);
			\draw(114.16,43.08) node[anchor=base west]{\fontsize{12.23}{15.07}\selectfont $G$};
			\draw(103.08,54.17) node[anchor=base west]{\fontsize{12.23}{15.07}\selectfont $H$};
			\draw(67.92,79.19) node[anchor=base west]{\fontsize{10.23}{13.07}\selectfont $V_{1H_0}$};
			\draw(69.07,65.82) node[anchor=base west]{\fontsize{10.23}{13.07}\selectfont $V_{2H_0}$};
			\draw(58.00,73.82) node[anchor=base west]{\fontsize{10.23}{13.07}\selectfont $H_0$};
\end{tikzpicture}%
 \caption{ The graph $G$ and $H_0\subseteq H \subseteq G$}\label{fig0}
\end{figure}
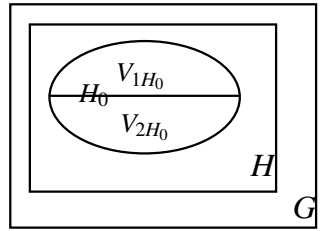
\begin{figure}
		\centering
		\begin{tikzpicture}[x=1.00mm, y=1.00mm, inner xsep=0pt, inner ysep=0pt, outer xsep=0pt, outer ysep=0pt]
			
		\end{tikzpicture}%
	\end{figure}

By the same proof as Claims \ref{claim0.1}, \ref{claim0.2}  and  \ref{findcycle},  we immediately obtain the following Claims \ref{claim0.3}, \ref{claim0.4}  and  \ref{claim0.5}, respectively.
\begin{claim}\label{claim0.3}
	For any  vertex $v_1, v_2\in V_{1H_0}$ (or $ v_1, v_2\in V_{2H_0}$), $|N_{H_0}(v_1)\cap N_{H_0}(v_2)|>k+5.$
\end{claim}
\begin{claim}\label{claim0.4}
Let $T\subseteq V(H_0)$ and $|T|\le 5.$	For any vertex $v_1\in V_{1H_0}$ and $ v_2\in V_{2H_0}$,  there exists a path $P_{v_1v_2}$ of order $2l$ for each integer $ l\in [2, k]$ in $H_0-T$.
\end{claim}
\begin{claim}\label{claim0.5}
For any vertex $u_1\in V_{1H_0}$ and $ v_1\in V_{2H_0}$, let $P_{u_1v_1}\subseteq G.$
 If $|V(P_{u_1v_1})|=5$ , then $G$ contains all odd cycles $C_{2l+3}$ for each integer $l\in[2,k-1]$.
\end{claim}

\begin{claim}\label{claim0.5'}
For any vertex $u_1, ~u_2\in V_{1H_0}$ (resp. $u_1, ~u_2\in V_{2H_0}$), let $P_{u_1u_2}\subseteq G.$
If $|V(P_{u_1u_2})|=6,$ then $G$ contains all odd cycles $C_{2l+3}$ for each integer $l\in[2,k-1]$.
 \end{claim}
\renewcommand\proofname{\bf Proof}
\begin{proof}
Without loss of generality, we only consider the situation that $u_1, ~u_2\in V_{1H_0}$.
Denote by  $P_{u_1u_2}:=u_1w_1w_2w_3w_4u_2$. By Claim \ref{claim0.3}, there exists a vertex $v_{21}\in V_{2H_0}$ such that $v_{21}\in N_{H_0}(u_1)\cap N_{H_0}(u_2)$ and $v_{21} \notin V(P_{u_1u_2}).$  Thus, we find a cycle $C_7:=u_2v_{21}u_1P_{u_1u_2}.$ Set $T=V(H_0)\cap (V(C_7)\setminus \{u_1, v_{21}\})$. Furthermore, by Claim \ref{claim0.4},  there exists a path $P'_{u_1v_{21}}$ of order $2l$ for each integer $ l\in [2, k-2]$ in $H_0-T$.
Therefore, we can find a $C_7$ and all odd cycles $C_{2l+5}:=v_{21}u_2w_4w_3w_2w_1u_1P'_{u_1v_{21}}$  for each integer $ l\in [2, k-2]$ in $G$, as desired.
\end{proof}

Since $\rho(G)\geq \rho(K_{\lceil \frac{n-2}{2} \rceil,\lfloor \frac{n-2}{2} \rfloor} \bullet K_3 )$ and $G\ncong K_{\lceil \frac{n-2}{2} \rceil,\lfloor \frac{n-2}{2} \rfloor} \bullet K_3 $, by Lemma \ref{nonbipartitec5}, $G$ contains  a cycle $C_5:=z_1z_2z_3z_4z_5z_1$  as a subgraph. In what follows, we shall prove that $G$ also contains all odd cycles  $C_{2l+3}$ for each integer $l\in[2,k-1]$ to complete the proof of Case 2.

Note that $|V(H_0)|>(1-40\theta)(n-4)$, $0< \theta\le \frac{1}{321}$ and $H_0$ is an induced subgraph of $G$. We deduce that
\begin{equation}\label{G-H0}
  |V(G-H_0)|=|V(G)|- |V(H_0)|< n-(1-40\theta)(n-4)<\frac{1}{8}n-10.
\end{equation}
Since $H_0$ is a bipartite graph, $H_0$ contains no odd cycle. Thus, $|V(C_5)\cap V(H_0)|\leq 4.$
Now we shall distinguish four situations as follows.
\begin{center}
  $ (II_1)~ 3\leq |V(C_5)\cap V(H_0)|\leq 4, ~(II_2)~|V(C_5)\cap V(H_0)|=2,$\\
\vspace{3mm}
  $(II_3)~|V(C_5)\cap V(H_0)|=1,~~~~~~(II_4)~|V(C_5)\cap V(H_0)|=0.$
\end{center}

For $(II_1)$ and $(II_2)$,  we can also find all odd cycles $C_{2l+3}$ for each integer $l\in [2,k-1]$ in $ G$ by Claim \ref{claim0.4}, Claim \ref{claim0.5} and using a similar proof as the Case 1 (i) and (ii), respectively.

For $(II_3)$, without loss of generality, let $V(C_5)\cap V(H_0)=V(C_5)\cap V_{1H_0}=\{z_5\}.$
Furthermore, we shall divide the proof of $(II_3)$  into the following two situations:
$$ (a_1) ~|V(C_5)\cap V(G- H)|=4, ~~(a_2)~0\leq |V(C_5)\cap V(G- H)|\leq 3.$$

For ($a_1$), note that $ |V(G-H)|=t\le 4.$ Then we see $t=4 $ and $V(G-H)=\{z_1, z_2, z_3, z_4\}$.
By (\ref{upperbound}) and  Lemma \ref{claim0.2'},
we have $d_G(z_i)\ge \frac{n}{8}-8$ for each $i\in\{1,2,3,4\}$.
Combining this with (\ref{G-H0}), for any $z\in \{z_1, z_2, z_3, z_4\}$,  we have
$$d_{H_0}(z)\geq d_G(z)-|V(G-H_0)|>\frac{n}{8}-8-(\frac{n}{8}-10)=2.$$
Therefore, by Claim \ref{claim0.3}, Claim \ref{claim0.5} and using a similar analysis as Case 1 (iii),  we complete the proof of ($a_1$).

For ($a_2$), it follows  $1\leq |V(C_5)\cap V(H-H_0)|\leq 4,$  and so there exists a vertex $z\in V(C_5)\cap V(H-H_0)$, which is adjacent to $z_5$ in $C_5$ or not, say $z=z_4$ or $z=z_3$. Since $d_H(z)\ge \delta(H)>\frac{n}{8}$ and $|V(H- H_0)|<|V(G- H_0)|<\frac{1}{8}n-10$ by (\ref{G-H0}), we deduce that
$$d_{H_0}(z)\geq d_H(z)-|V(H- H_0)|>\frac{n}{8}-(\frac{1}{8}n-10)>2.$$
Thus, if $z=z_4$, there exists a vertex $z_0\in N_{H_0}(z_4)\setminus \{z_5\}$.
 When $z_0\in  V_{1H_0},$ then we can find a path $P_{z_5z_0}:=z_5z_1z_2z_3z_4z_0$ of order $6$ satisfying $z_0,z_5\in V_{1H_0}.$ By Claim \ref{claim0.5'}, the result follows.
When $z_0\in  V_{2H_0},$ by Claim \ref{claim0.4}, we can find a path $P_{z_5z_0}\subseteq H_0$ of order $2l$ for each integer $ l \in [2,k]$, and so $G$ contains all odd cycles $C_{2l+1}:=z_0z_4z_5P_{z_5z_0}$ for each $l\in [2,k] $, as required.
If $z=z_3$, then we can use a similar analysis to complete the proof.

For $(II_4)$,
there exists at least one vertex belonging to $ V(C_5)\cap V(H-H_0)$ since $t=n-h\le 4$.
In what follows, we shall divide the proof into the following three situations:
  $$(b_1)~ |V(C_5)\cap V(H-H_0)|=1, ~(b_2)~|V(C_5)\cap V(H-H_0)|\geq 3,
  ~(b_3)~|V(C_5)\cap V(H- H_0)|=2.$$

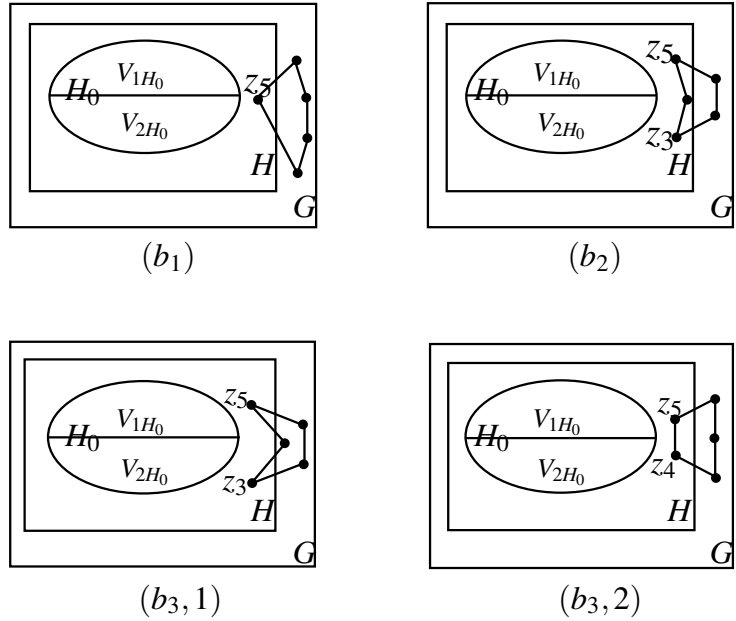
\begin{figure}[hpt]
\centering
\begin{tikzpicture}[x=1.00mm, y=1.00mm, inner xsep=0pt, inner ysep=0pt, outer xsep=0pt, outer ysep=0pt,scale=0.5]
\path[line width=0mm] (70.81,-40.83) rectangle +(194.09,160.65);
		\definecolor{L}{rgb}{0,0,0}
		\path[line width=0.30mm, draw=L] (39.84,40.41) rectangle +(80.44,59.23);
		\path[line width=0.30mm, draw=L] (44.99,49.96) rectangle +(64.77,44.33);
		\path[line width=0.30mm, draw=L] (75.18,74.99) ellipse (25.03mm and 14.90mm);
		\path[line width=0.30mm, draw=L] (49.96,75.37) -- (100.21,75.37);
			\draw(114.16,43.08) node[anchor=base west]{\fontsize{12.23}{15.07}\selectfont $G$};
			\draw(103.08,54.17) node[anchor=base west]{\fontsize{12.23}{15.07}\selectfont $H$};
			\draw(67.92,79.19) node[anchor=base west]{\fontsize{10.23}{13.07}\selectfont $V_{1H_0}$};
			\draw(69.07,65.82) node[anchor=base west]{\fontsize{10.23}{13.07}\selectfont $V_{2H_0}$};
			\path[line width=0.30mm, draw=L] (149.68,40.41) [rotate around={360:(149.68,40.41)}] rectangle +(80.22,59.41);
			\path[line width=0.30mm, draw=L] (154.62,49.96) rectangle +(64.77,44.33);
			\path[line width=0.30mm, draw=L] (184.81,74.99) ellipse (25.03mm and 14.90mm);
			\path[line width=0.30mm, draw=L] (159.59,75.37) -- (209.84,75.37);
			\draw(223.79,43.08) node[anchor=base west]{\fontsize{12.23}{15.07}\selectfont $G$};
			\draw(212.71,54.17) node[anchor=base west]{\fontsize{12.23}{15.07}\selectfont $H$};
			\draw(177.55,79.19) node[anchor=base west]{\fontsize{10.23}{13.07}\selectfont $V_{1H_0}$};
			\draw(178.70,65.82) node[anchor=base west]{\fontsize{10.23}{13.07}\selectfont $V_{2H_0}$};
			\path[line width=0.30mm, draw=L] (39.81,-49.57) [rotate around={360:(39.81,-49.57)}] rectangle +(80.13,59.60);
			\path[line width=0.30mm, draw=L] (43.60,-40.10) rectangle +(65.96,45.45);
			\path[line width=0.30mm, draw=L] (74.78,-15.30) ellipse (25.03mm and 14.90mm);
			\path[line width=0.30mm, draw=L] (49.96,-15.32) -- (100.21,-15.32);
			\draw(114.16,-48.61) node[anchor=base west]{\fontsize{12.23}{15.07}\selectfont $G$};
			\draw(103.08,-37.52) node[anchor=base west]{\fontsize{12.23}{15.07}\selectfont $H$};
			\draw(67.92,-12.50) node[anchor=base west]{\fontsize{10.23}{13.07}\selectfont $V_{1H_0}$};
			\draw(69.07,-25.87) node[anchor=base west]{\fontsize{10.23}{13.07}\selectfont $V_{2H_0}$};
			\path[line width=0.30mm, draw=L] (150.24,-49.97) rectangle +(79.66,59.40);
			\path[line width=0.30mm, draw=L] (155.02,-39.93) rectangle +(64.77,44.33);
			\path[line width=0.30mm, draw=L] (184.61,-15.10) [rotate around={0:(184.61,-15.10)}] ellipse (25.03mm and 14.90mm);
			\path[line width=0.30mm, draw=L] (159.39,-15.52) -- (209.64,-15.52);
			\draw(223.79,-48.61) node[anchor=base west]{\fontsize{12.23}{15.07}\selectfont $G$};
			\draw(212.71,-37.52) node[anchor=base west]{\fontsize{12.23}{15.07}\selectfont $H$};
			\draw(177.55,-12.50) node[anchor=base west]{\fontsize{10.23}{13.07}\selectfont $V_{1H_0}$};
			\draw(178.70,-25.87) node[anchor=base west]{\fontsize{10.23}{13.07}\selectfont $V_{2H_0}$};
			\definecolor{F}{rgb}{0,0,0}
			\path[line width=0.30mm, draw=L, fill=F] (115.14,84.63) circle (1.00mm);
			\path[line width=0.30mm, draw=L, fill=F] (117.68,74.76) circle (1.00mm);
			\path[line width=0.30mm, draw=L, fill=F] (117.96,64.05) circle (1.00mm);
			\path[line width=0.30mm, draw=L, fill=F] (115.43,54.75) circle (1.00mm);
			\path[line width=0.30mm, draw=L, fill=F] (105.00,74.20) circle (1.00mm);
			\path[line width=0.30mm, draw=L] (115.14,84.63) -- (105.00,75.04);
			\path[line width=0.30mm, draw=L] (114.58,84.91) -- (117.68,75.04);
			\path[line width=0.30mm, draw=L] (117.93,74.53) -- (117.95,64.33);
			\path[line width=0.30mm, draw=L] (118.25,64.33) -- (115.43,54.75);
			\path[line width=0.30mm, draw=L] (104.71,75.04) -- (115.53,54.40);
			\draw(101.26,76.61) node[anchor=base west]{\fontsize{14.23}{17.07}\selectfont $z_5$};
			\path[line width=0.30mm, draw=L, fill=F] (214.88,84.97) circle (1.00mm);
			\path[line width=0.30mm, draw=L, fill=F] (217.87,74.21) circle (1.00mm);
			\path[line width=0.30mm, draw=L, fill=F] (215.08,64.24) circle (1.00mm);
			\path[line width=0.30mm, draw=L, fill=F] (225.44,79.79) circle (1.00mm);
			\path[line width=0.30mm, draw=L, fill=F] (225.24,70.02) circle (1.00mm);
			\path[line width=0.30mm, draw=L] (215.08,84.97) -- (225.24,79.59);
			\path[line width=0.30mm, draw=L] (225.64,79.99) -- (225.64,69.82);
			\path[line width=0.30mm, draw=L] (214.88,64.24) -- (225.44,69.82);
			\path[line width=0.30mm, draw=L] (214.68,84.97) -- (217.87,74.41);
			\path[line width=0.30mm, draw=L] (217.67,74.41) -- (215.08,64.44);
			\draw(207.00,62.45) node[anchor=base west]{\fontsize{14.23}{17.07}\selectfont $z_3$};
			\draw(207.00,85.58) node[anchor=base west]{\fontsize{14.23}{17.07}\selectfont $z_5$};
			\path[line width=0.30mm, draw=L, fill=F] (103.25,-6.72) circle (1.00mm);
			\path[line width=0.30mm, draw=L, fill=F] (112.02,-16.88) circle (1.00mm);
			\path[line width=0.30mm, draw=L, fill=F] (103.45,-27.45) circle (1.00mm);
			\path[line width=0.30mm, draw=L, fill=F] (116.81,-11.90) circle (1.00mm);
			\path[line width=0.30mm, draw=L, fill=F] (117.01,-22.46) circle (1.00mm);
			\path[line width=0.30mm, draw=L] (103.45,-6.72) -- (116.41,-11.90);
			\path[line width=0.30mm, draw=L] (117.21,-11.70) -- (117.21,-21.87);
			\path[line width=0.30mm, draw=L] (104.05,-27.25) -- (117.41,-22.06);
			\path[line width=0.30mm, draw=L] (103.05,-6.72) -- (112.22,-16.88);
			\path[line width=0.30mm, draw=L] (112.02,-17.08) -- (103.45,-27.25);
			\draw(95.68,-29.24) node[anchor=base west]{\fontsize{12.23}{15.07}\selectfont $z_3$};
			\draw(96.07,-6.11) node[anchor=base west]{\fontsize{12.23}{15.07}\selectfont $z_5$};
			\draw(54.22,74.06) node[anchor=base west]{\fontsize{12.23}{15.07}\selectfont $H_0$};
			\draw(161.86,74.06) node[anchor=base west]{\fontsize{12.23}{15.07}\selectfont $H_0$};
			\draw(54.22,-17.63) node[anchor=base west]{\fontsize{12.23}{15.07}\selectfont $H_0$};
			\draw(161.86,-17.63) node[anchor=base west]{\fontsize{12.23}{15.07}\selectfont $H_0$};
			\path[line width=0.30mm, draw=L, fill=F] (214.68,-10.55) circle (1.00mm);
			\path[line width=0.30mm, draw=L, fill=F] (214.88,-20.12) circle (1.00mm);
			\path[line width=0.30mm, draw=L, fill=F] (225.24,-5.17) circle (1.00mm);
			\path[line width=0.30mm, draw=L, fill=F] (225.04,-15.54) circle (1.00mm);
			\path[line width=0.30mm, draw=L, fill=F] (225.44,-26.10) circle (1.00mm);
			\path[line width=0.30mm, draw=L] (215.08,-10.35) -- (225.44,-4.97);
			\path[line width=0.30mm, draw=L] (214.68,-10.35) -- (214.68,-19.92);
			\path[line width=0.30mm, draw=L] (225.24,-5.37) -- (225.24,-15.74);
			\path[line width=0.30mm, draw=L] (225.24,-15.54) -- (225.24,-25.90);
			\path[line width=0.30mm, draw=L] (214.28,-19.92) -- (224.84,-25.50);
			\draw(210.09,-8.56) node[anchor=base west]{\fontsize{12.23}{15.07}\selectfont $z_5$};
			\draw(208.30,-24.71) node[anchor=base west]{\fontsize{12.23}{15.07}\selectfont $z_4$};
			\draw(74.15,30.09) node[anchor=base west]{\fontsize{12.23}{17.07}\selectfont $(b_1)$};
			\draw(186.37,30.09) node[anchor=base west]{\fontsize{12.23}{17.07}\selectfont $(b_2)$};
			\draw(73.55,-60.75) node[anchor=base west]{\fontsize{12.23}{17.07}\selectfont $(b_3,1)$};
			\draw(183.98,-60.75) node[anchor=base west]{\fontsize{12.23}{17.07}\selectfont $(b_3,2)$};
\end{tikzpicture}%
\caption{ some situations of ($II_4$)}\label{fig1}
\end{figure}
For ($b_1$), without loss of generality, suppose $\{z_5\}= V(C_5)\cap V(H-H_0)$, which implies that $t=4 $ and $V(G-H)=\{z_1, z_2, z_3, z_4\}$(\emph{see Figure \ref{fig1}:($b_1$)}).
 By (\ref{upperbound}) and Lemma \ref{claim0.2'},
we have $d_G(z_i)\ge \frac{n}{8}-8$ for each $i\in\{1,2,3,4\}$.
Recall (\ref{G-H0}) that  $|V(G- H_0)|<\frac{1}{8}n-10$. Then,
$$d_{H_0}(z_3)\ge d_G(z_3)-|V(G- H_0)|>2. $$
Thus, there exists a vertex $z_0\in N_{H_0}(z_3)$. 
Without loss of generality, we suppose $z_0\in  V_{1H_0}$.
Also, since $z_5 \in V(H-H_0)$ and $\delta(H)>\frac{n}{8}$, we have
$$d_{H_0}(z_5)\geq d_G(z_5)-|V(G-H_0)|\geq \frac{n}{8}- (\frac{1}{8}n-10)>2,$$
and so there exists some vertex $v \in N_{H_0}(z_5)\setminus \{z_0\}$.
If $v \in V_{1H_0}$, then  we can find a path $P_{vz_0}:=vz_5z_1z_2z_3z_0$ of order $6$ satisfying $z_0,~v\in V_{1H_0}.$ By Claim \ref{claim0.5'}, the result follows.
If $v \in  V_{2H_0},$ then we can find a path $P_{vz_0}:=vz_5z_4z_3z_0$ of order $5$ satisfying $z_0\in V_{1H_0}$ and $v\in V_{2H_0}.$
By Claim \ref{claim0.5},  we complete the proof of ($b_1$).

For ($b_2$), there exist two non-adjacent vertices in $C_5$, say $\{z_3,  z_5\}$, belonging to $V(C_5)\cap V(H-H_0)$ (\emph{see Figure \ref{fig1}:($b_2$)}).
Note that $\delta(H)>\frac{n}{8}$ and $|V(H-H_0)|\le |V(G-H_0)|<\frac{1}{8}n-10$ by (\ref{G-H0}).
Then we deduce that for each $i\in\{3,5\}$,
$$ d_{H_0}(z_i)\ge d_H(z_i)-|V(H-H_0)|\ge \delta(H)-|V(H-H_0)|>2.  $$
 Without loss of generality, there exists some vertex $v_{11}\in N_{H_0}(z_5)\cap V_{1H_0}$, and let $v\in N_{H_0}(z_3) \setminus \{v_{11}\}.$
  When $v\in V_{1H_0}$, then  we can find a path $P_{vv_{11}}:=vz_3z_2z_1z_5v_{11}$ of order $6$ satisfying $v_{11},v\in V_{1H_0}.$ By Claim \ref{claim0.5'}, the result follows.
 When $v \in V_{2H_0},$ then we can find a path $P_{v_{11}v}:=vz_3z_4z_5v_{11}$ of order $5$ satisfying $v_{11}\in V_{1H_0}$ and $v\in V_{2H_0}$. By Claim \ref{claim0.5},  we complete the proof of ($b_2$).

For ($b_3$), then $t=3$ or $4.$
 If $t=4$, then by (\ref{upperbound}) and Lemma \ref{claim0.2'}, we see $d_G(v)\ge \frac{n}{8}-8$ for any $v\in V(G-H)$. Recall (\ref{G-H0}). Then  any  vertex $v\in V(G-H)$ satisfies
 $$d_{H_0}(v)> d_G(v)-|V(G-H_0)|>2.$$
Note that the condition $|V(C_5)\cap V(H_0)|=0$ and $|V(C_5)\cap V(H-H_0)|=2$. Then $ |V(C_5)\cap V(G-H)|=3$, and so we know that there exist two non-adjacent vertices in $C_5$, say $\{z_3,  z_5\}$, belonging to $V(C_5)\cap V(G-H)$. Thus, by using a similar analysis  as ($b_2$), the result follows.

Now, we consider that $t=3$. Without loss of generality, we consider $\{z_5, z_i\}= V(C_5)\cap V(H-H_0)$ for $i\in \{3,4\}$. Note that $\delta(H)>\frac{n}{8}$ and $|V(G-H_0)|<\frac{1}{8}n-10$ by (\ref{G-H0}). When $z_i=z_3$, (\emph{see Figure \ref{fig1}:($b_3$.1)}), then by using a similar analysis as ($b_2$), we obtain all odd cycles $C_{2l+1}$ for each integer $l\in [2,k]$ in $G$, as required.
When $z_i=z_4$ (\emph{see Figure \ref{fig1}:($b_3.2$)}), then $\{z_1, z_2, z_3\}=V(G-H)$ and $z_4, z_5 \in V(H-H_0)$.   Note that $\delta(H)>\frac{n}{8}$ and  $|V(H- H_0)|\le |V(G- H_0)|<\frac{1}{8}n-10$ due to (\ref{G-H0}). 
Then we deduce that for each $i\in\{4,5\}$,
$$ d_{H_0}(z_i)\ge d_H(z_i)-|V(H-H_0)|\ge \delta(H)-|V(H-H_0)|>2.  $$
Therefore, we can  distinguish the following two situations: there exist two different vertices $u$, $v$  such that
\begin{itemize}
  \item[($\ast$)]  $u\in N_{H_0}(z_4)\cap V_{1H_0}$ and $v\in N_{H_0}(z_5)\cap V_{1H_0}$,
  or
  \item[($\diamond$)] $u\in N_{H_0}(z_4)\cap V_{1H_0}$ and $v\in N_{H_0}(z_5)\cap V_{2H_0}$.
\end{itemize}

 For ($\ast$),  by Claim \ref{claim0.3},  there exists a vertex $v_{21}\in N_{H_0}(u)\cap N_{H_0}(v)$.
 Then we can find a path $P_{v_{21}v}:=vz_5z_4uv_{21}$ of order $5$ satisfying $v_{21}\in V_{2H_0}$ and $v\in V_{1H_0}$. By Claim \ref{claim0.5},  the result follows.

For ($\diamond$),
we distinguish the degree of $z_2$ to complete our proof.
If $d_G(z_2)\le 4$, since $d_{G-\{z_i|1\le i\le j\}}(z_{j+1})\le d_G (z_{j+1})$ for each $j\in\{1,2\}$, by Lemma \ref{deletingvertex0}, we have
\begin{align*}
   &\rho^2(G-\{z_1\}) \geq \rho^2(G) - 2d_G(z_1),\\
  &\rho^2(G-\{{z_1}, {z_2}\}) \geq \rho^2(G-\{z_1\})-2d_G(z_2)\ge \rho^2(G-\{z_1\})-8, \\
  &\rho^2(H)=\rho^2(G-\{{z_1}, {z_2},  {z_3}\}) \geq \rho^2(G-\{{z_1}, {z_2}\})-2d_G(z_3),
\end{align*}
which gives
$8+2(d_G(z_1)+d_G(z_3))\geq \rho^2(G)-\rho^2(H).$
Note that $t=3$.  By $h=n-3$, (\ref{lowerboundG}) and (\ref{upperbound}), we deduce that
$$d_G(z_1)+d_G(z_3)\ge \frac{1}{2}( \rho^2(G)-\rho^2(H)-8)\geq \frac{n^2-4n+3}{8}-\frac{(n-3)^2}{8}-4=\frac{2n-38}{8},$$
which implies that  $\max \{d_G(z_1), d_G(z_3)\}>\frac{n}{8}-8$. Without loss of generality, set $d_G(z_1)>\frac{n}{8}-8.$
 Combining this with  (\ref{G-H0}), we have
  $$d_{H_0}(z_1)\geq d_G(z_1)-|V(G-H_0)|> \frac{n}{8}-8-(\frac{n}{8}-10)=2.$$
Hence, there exists a vertex $w\in N_{H_0}(z_1)\setminus \{u, v\}.$
If $ w\in V_{1H_0}$,
then  we can find a path $P_{uw}:=uz_4z_3z_2z_1w$ of order $6$ satisfying  $u,~w\in V_{1H_0}.$ By Claim \ref{claim0.5'},
the result follows.
If $ w\in V_{2H_0}$, then  we can find a path $P_{uw}:=uz_4z_5z_1w$ of order $5$ satisfying $u\in V_{1H_0}$ and $w\in V_{2H_0}$. By Claim \ref{claim0.5}, the result follows.

 If $d_G(z_2)\geq 5$, then there exists a vertex $z_0\in N_H(z_2)\setminus\{z_4, z_5\}$.
Furthermore, if $z_0\in V(H-H_0)$, then
$$d_{H_0}(z_0)\geq d_H(z_0)-|V(H-H_0)|> \frac{n}{8}-(\frac{n}{8}-10)>2$$
since $\delta(H)>\frac{n}{8}$ and
$|V(H-H_0)|\leq|V(G-H_0)|< \frac{n}{8}-10$ by (\ref{G-H0}).
Hence, there exists a vertex $u_1\in N_{H_0}(z_0)\setminus \{u\}.$
If $u_1\in  V_{1H_0},$
then  we can find a path $P_{uu_1}:=uz_4z_3z_2z_0u_1$ of order $6$ satisfying  $u,~u_1\in V_{1H_0}.$ By Claim \ref{claim0.5'}, the result follows.
If $u_1\in  V_{2H_0},$ by using a similar analysis, the result follows.
If $z_0\in V_{1H_0}$, then we can find a path $P_{uz_0}:=uz_4z_5z_1z_2z_0$ of order $6$ satisfying  $u,~z_0\in V_{1H_0}.$ By Claim \ref{claim0.5'}, the result follows. If $z_0\in V_{2H_0}$, by using a similar analysis, we complete the proof of ($b_3$).

To sum up, we complete the proof.


\begin{thebibliography}{99}
\bibitem{BG}  L. Babai, B. Guiduli, Spectral extrema for graphs: the Zarankiewicz problem, \emph{Electron. J. Combin.,} \textbf{16 } (2009) \#R123.

\bibitem{CDT22+}
S. Cioab\u{a}, D.N. Desai, M. Tait, The spectral even cycle problem, arXiv:2205.00990.


\bibitem{GH19}
J. Gao, X. Hou, The spectral radius of graphs without long cycles, \emph{Linear Algebra Appl.} \textbf{566 } (2019) 17--33.

\bibitem{GN20}
J. Ge, B. Ning, Spectral radius and Hamiltonian properties of graphs, II, \emph{Linear Multilinear Algebra}  \textbf{68}  (2020) 2298--2315.


\bibitem{GWL2021+}
J. Guo, Z. Wang, Z. Lou, X. Li, Improved upper bound for the spectral radius of a graph. (submitted)

\bibitem{GLZ}
H. Guo, H. Lin, Y. Zhao, A spectral condition for the existence of a pentagon in non-bipartite graphs, \emph{Linear Algebra Appl.}  \textbf{627} (2021) 140--149.

\bibitem{LN}
B. Li, B. Ning, Eigenvalues and cycles of consecutive lengths, arXiv:2110.05670.

\bibitem{LNW21+}
H. Lin, B. Ning, B. Wu, Eigenvalues and triangles in graphs, \emph{Combin. Probab. Comput.}  \textbf{30}  (2021) 258--270.

\bibitem{LF22+}
Y. Li, W. Liu, L. Feng, A survey on spectral conditions for some extremal graph problems, arXiv:2111.03309.


\bibitem{NP20}
B. Ning, X. Peng, Extensions of the Erd\H{o}s-Gallai theorem and Luo's theorem, \emph{Combin. Probab. Comput.} \textbf{29} (2020), no. 1, 128--136.

\bibitem{NP2006}
V. Nikiforov, R.H. Schelp, Cycle lengths in graphs with large minimum degree, \emph{J. Graph Theory} \textbf{52} (2006) 157--170.

\bibitem{N07}  V. Nikiforov, Bounds on graph eigenvalues II, \emph{Linear Algebra Appl.} \textbf{427} (2007) 183--189.

\bibitem{N08}
V. Nikiforov, A spectral condition for odd cycles in graphs, \emph{Linear Algebra Appl.} \textbf{428} (2008) 1492--1498.


\bibitem{N09}
V. Nikiforov, The maximum spectral radius of $C_4$-free graphs of given order and size,
\emph{Linear Algebra Appl.} \textbf{430} (2009) 2898--2905.

\bibitem{NI4}  V. Nikiforov, A contribution to the Zarankiewicz problem, \emph{Linear Algebra Appl.,} \textbf{432} (2010)  1405--1411.

\bibitem{N10}
V. Nikiforov, The spectral radius of graphs without paths and cycles of specified length, \emph{Linear Algebra Appl.} \textbf{432} (2010) 2243--2256.

\bibitem{N11}
V. Nikiforov, Some new results in extremal graph theory. Surveys in combinatorics 2011, 141--181,
London Math. Soc. Lecture Note Ser., 392, Cambridge Univ. Press, Cambridge, 2011.



\bibitem{WI} H. Wilf, Spectral bounds for the clique and independence numbers of graphs, \emph{J. Combin. Theory Ser. B.} \textbf{40} (1986) 113--117.


\bibitem{ZL20}
M. Zhai, H. Lin, Spectral extrema of graphs: forbidden hexagon, \emph{Discrete Math.} \textbf{343} (2020), no. 10, 112028, 6 pp.


\bibitem{ZLar}
M. Zhai,  H. Lin, A strengthening of the spectral color critical edge theorem: books and theta graphs, 2021,
 arXiv:2102.04041.

\bibitem{ZW12}
M. Zhai, B. Wang, Proof of a conjecture on the spectral radius of $C_4$-free graphs,
\emph{Linear Algebra Appl.} \textbf{437} (2012) 1641--1647.




\end{thebibliography}
\end{document}